\documentclass{article}%
\usepackage{amsmath}
\usepackage{amsfonts}
\usepackage{amssymb}
\usepackage{graphicx}%
\providecommand{\U}[1]{\protect\rule{.1in}{.1in}}

\newenvironment{proof}[1][Proof]{\noindent\textbf{#1.} }{\ \rule{0.5em}{0.5em}}
\begin{document}
\ {\Large Compact Group Actions On Operator Algebras}

{\Large \ \ \ \ \ \ \ \ \ \ \ \ \ \ \ \ \ \ \ \ and Their Spectra}

\bigskip

{\large Costel Peligrad}

\bigskip

\textit{Department of Mathematical Sciences, University of Cincinnati, 4508
French Hall West, Cincinnati, OH 45221-0025, United States. E-mail address:
costel.peligrad@uc.edu}

\bigskip

{\large Abstract. }We consider a class of dynamical systems with compact non
abelian groups that include C*-, W*- and multiplier dynamical systems. We
prove results that relate the algebraic properties such as simplicity or
primeness of the fixed point algebras as defined in Section 3., to the
spectral properties of the action, including the Connes and strong Connes spectra.

\bigskip


\bigskip


\section{\bigskip\bigskip Introduction}

In [2], Connes introduced the invariant $\Gamma(U)$ known as the Connes
spectrum of the action $U$\ of a locally compact abelian group on a von
Neumann algebra\ and used it in his seminal classification of type III von
Neumann factors. Soon after, Olesen [10] defined the Connes spectrum of an
action of a locally compact abelian group on a C*-algebra. In [11], using the
definition of the Connes spectrum in [10], it is proven an analog of a result
of Connes and Takesaki [3, Chapter III, Corollary 3.4.] regarding the
significance of the Connes spectrum of a locally compact abelian group action
on a C*-algebra for the ideal structure of the crossed product. In particular,
in [11] is discussed a spectral characterization for the crossed product to be
a prime C*-algebra. This definition of the Connes spectrum in [10] cannot be
used to prove similar results for the simplicity of the crossed product,
unless the group is discrete [11]. Kishimoto [8] defined the strong Connes
spectrum for C*-dynamical systems with locally compact abelian groups that
coincides with the Connes spectrum for the W*-dynamical system and with the
Connes spectrum defined by Olesen for discrete abelian group actions on
C*-algebras and he proved the Connes-Takesaki result for simple crossed
products. In [2] Connes obtained results that relate the spectral properties
of the von Neumann algebra with the algebraic properties of the fixed point
algebra. These results were extended in [12] to C*-algebras and compact
abelian groups. In [6], [14] we considered the problems of simplicity and
primeness of the crossed product by compact, non abelian group actions. In
particular, in [6] we have defined the Connes and strong Connes spectra for
such actions that coincide with Connes spectra [2], [10], respectively with
the strong Connes spectra [8] for compact abelian groups. Further, in [15] we
have considered the case of one-parameter $\mathcal{F}$-dynamical systems that
include the C*- the W*- and the multiplier one-parameter dynamical systems. In
particular, we have obtained extensions of some results in [2], [12] for
$W^{\ast}-,$ respectively $C^{\ast}-$ dynamical systems to the case
of$\ \mathcal{F}$-dynamical systems with compact abelian groups [15, Theorems
3.2 and 3.4.]. In this paper we will prove results for $\mathcal{F}$-dynamical
systems with compact non abelian groups. Our results contain and extend to the
case of compact non abelian groups\ the following: [2 Proposition 2.2.2. b)
and Theorem 2.4.1], [12, Theorem 2], [13, Theorem 8.10.4] and [15, Theorems
3.2. and 3.4.]. In Section 2. we will set up the framework and state some
results that will be used in the rest of the paper. In Section 3. we discuss
the connection between the strong Connes spectrum, $\widetilde{\Gamma
}_{\mathcal{F}}(\alpha),$ of the action and the $\mathcal{F}$-simplicity of
the fixed point algebras $(X\otimes B(H_{\pi}))^{\alpha\otimes ad\pi}$. In
Section 4. we will get similar results about the connection between the
$\mathcal{F}$-primeness of the fixed point algebras and the Connes spectrum,
$\Gamma_{\mathcal{F}}(\alpha),$ of the action.

\section{\bigskip Notations and preliminary results}

\bigskip

This section contains the definitions of the basic concepts used in the rest
of the paper, the notations and some preliminary results.

\bigskip

\textbf{2.1. Definition. ([1], [16 ]) }\textit{A dual pair of Banach spaces
is, by definition, a pair} $(X,\mathcal{F})$ \textit{of Banach spaces with the
following properties:}

\textit{a)} $\mathcal{F}$ \textit{is a Banach subspace of the dual} $X^{\ast}%
$\ \textit{of} $X.$

\textit{b)} $\left\Vert x\right\Vert =\sup\left\{  \left\vert \varphi
(x)\right\vert :\varphi\in\mathcal{F},\left\Vert \varphi\right\Vert
\leq1\right\}  ,x\in X.$

\textit{c)} $\left\Vert \varphi\right\Vert =\sup\left\{  \left\vert
\varphi(x)\right\vert :x\in X,\left\Vert x\right\Vert \leq1\right\}
,\varphi\in\mathcal{F}.$

\textit{d)} \textit{The convex hull of every relativel}y $\mathcal{F}%
$-\textit{compact subset of }$\mathit{X}$\textit{\ is relatively}
$\mathcal{F}$-\textit{compact}.

\textit{e)} \textit{The convex hull of every relatively }$X$\textit{-compact
subset of} $\mathcal{F}$\ \textit{is relatively} $X$-\textit{compact}.

\textit{In the rest of the paper }$X$\textit{\ will be assumed to be a
C*-algebra with the additional property}

\textit{f) The involution of }$X$\textit{\ is }$\mathcal{F}$%
\textit{-continuous and the multiplication in }$X$\textit{\ is separately
}$\mathcal{F}$\textit{-continuous.}

\bigskip

The property d) implies the existence of the weak integrals of continuous
functions defined on a locally compact measure space, ($S,\mu)$ with values in
$X$ endowed with the $\mathcal{F}$-topology:

If $f$ is such a function, we will denote by%
\[
\int f(s)d\mu
\]
the unique element $y$\ of $X$\ such that
\[
\varphi(y)=\int_{S}\varphi(f(s))d\mu
\]
for every $\varphi\in\mathcal{F}$ [1, Proposition 1.2.]$.$\ The propery e) was
used by Arveson [1, Proposition 1.4.] to prove the continuity in the
$\mathcal{F}$-topology of some linear mappings on $X\ $(in particular the
mappings $P_{\alpha}(\pi)$\ and ($P_{\alpha})_{ij}(\pi)$ defined below).

\bigskip

\textbf{2.2}. \textbf{Examples.} \textit{a) [1] If }$X$\textit{\ is a
C*-algebra and and }$\mathcal{F}=X^{\ast},$\textit{ conditions 1)-5) are
satisfied.}

\textit{b) [1] If }$X$\textit{\ is a W*-algebra and }$\mathcal{F}=X_{\ast}%
$\textit{ is its predual then conditions 1)-5) are satisfied.}

\textit{c) [4] If }$X=M(Y)$\textit{\ is the multiplier algebra of }%
$Y$\textit{\ and }$\mathcal{F}=Y^{\ast}$\textit{\ then conditions 1)-5) are
satisfied. In addition, in this case, the }$\mathcal{F}$\textit{-topology on
}$X$\textit{\ is compatible with the strict toplogy on }$X=M(Y).$

\bigskip

Let $(X,\mathcal{F})$ be a dual pair of Banach spaces $G$ a compact group and
$\alpha:G\rightarrow Aut(X)$ a homeomorphism of $G$ into the group of $\ast
-$automorphisms of $X$. We say that $(X,G,\alpha)$ is an $\mathcal{F}%
$-dynamical system if the mapping%
\[
g\rightarrow\varphi(\alpha_{g}(x))
\]
is continuous for every $x\in X$\ and $\varphi\in\mathcal{F}.$

\bigskip

\textbf{2.3}. \textbf{Examples. }\textit{a) If }$\mathcal{F}=X^{\ast},$ the
dual of $X$\textit{ then, by [7 p. 306] the above condition is equivalent to
the continuity of }the mapping $g\rightarrow\alpha_{g}(x)$\ from
$G$\textit{\ to }$X$\textit{\ endowed with the norm topology for every }$x\in
X$\textit{, so, in this case }$(X,G,\alpha)$\textit{ is a C*-dynamical
system.}

\textit{b) If }$X$\textit{\ is a von Neumann algebra and }$\mathcal{F}%
=X_{\ast},$ the predual of $X$\textit{ then }$(X,G,\alpha)$\textit{ is a
W*-dynamical system.}

\textit{c) If }$X=M(Y)$\textit{\ is the multiplier algebra of }$Y$%
\textit{\ and }$\mathcal{F}=Y^{\ast},$\textit{ then }$(X,G,\alpha)$\textit{ is
said to be a multiplier dynamical system.}

\bigskip

Let $(X,G,\alpha)$ be an $\mathcal{F}$-dynamical system with $G$\ compact.
Denote by $\widehat{G}$ the set of unitary equivalence classes of irreducible
representations of $G.$ For each $\pi\in\widehat{G}$\ denote also by $\pi$\ a
fixed representative of that class. If $\chi_{\pi}(g)=d_{\pi}\sum
_{i=1}^{d_{\pi}}\pi_{ii}(g^{-1})=d_{\pi}\sum\overline{\pi_{ii}(g)}$ is the
character of $\pi,$\ denote by
\[
P_{\alpha}(\pi)(x)=\int_{G}\chi_{\pi}(g)\alpha_{g}(x)dg.
\]
Then $P_{\alpha}(\pi)$ is a projection of $X$\ onto the spectral subspace
\[
X_{1}(\pi)=\left\{  x\in X:P_{\alpha}(\pi)(x)\right\}  .
\]
where the integral is taken in the weak sense\ defined in (1) above. As in
[14] one can also define for every $1\leq i,j\leq d_{\pi}$%
\[
(P_{\alpha})_{ij}(\pi)(x)=\int_{G}\overline{\pi_{ji}(g)}\alpha_{g}(x)dg.
\]
where\ $d_{\pi}$ is the dimension of the Hilbert space $H_{\pi}$\ of $\pi
$\ and show that
\[
(P_{\alpha})_{ij}(\pi)(X)\subset X_{1}(\pi).
\]
Using [1, Proposition 1.4.] it follows that $P_{\alpha}(\pi),(P_{\alpha}%
)_{ij}(\pi)$\ are $\mathcal{F}$-continuous. If $\pi$\ is the identity one
dimensional representation $\iota$\ of $X,$\ we will denote
\[
P_{\alpha}(\iota)=P_{\alpha}.
\]
and
\[
X_{1}(\iota)=X^{\alpha}.
\]
is the fixed point algebra of the action.

\bigskip

\textbf{2.4. Remark. } $\overline{\sum_{\pi\in\widehat{G}}X_{1}(\pi)}^{\sigma
}=X,$\textit{ where }$\overline{\sum_{\pi\in\widehat{G}}X_{1}(\pi)}^{\sigma}%
$\textit{\ denotes the closure of }$\sum_{\pi\in\widehat{G}}X_{1}(\pi
)$\textit{ in the }$\mathcal{F}$\textit{-topology of }$X.$

\bigskip

\begin{proof}
Suppose that there exists $\varphi\in\mathcal{F}$\ such that $\varphi
(X_{1}(\pi))=\left\{  0\right\}  $ for every $\pi\in\widehat{G}.$ Since, as
noticed above, $(P_{\alpha})_{ij}(\pi)(X)\subset X_{1}(\pi),$ it follows that%
\[
\int_{G}\overline{\pi_{ij}(g)}\varphi(\alpha_{g}(x))dg=0.
\]
for every $x\in X$ and every $\pi\in\widehat{G}.$ Since $\left\{  \pi
_{ij}(g):\pi\in\widehat{G},1\leq i,j\leq d_{\pi}\right\}  $ is an orthogonal
basis of $L^{2}(G),$ and $\varphi(\alpha_{g}(x))$\ is a continuous function of
$g,$\ for every $x\in X,$\ it follows that $\varphi(x)=0$ for every $x\in
X$\ so $\varphi=0$ and we are done.
\end{proof}

\bigskip

In ([9], [14], [6]) it is pointed out that the spectral subspaces%
\[
X_{2}(\pi)=\left\{  a\in X\otimes B(H_{\pi}):(\alpha_{g}\otimes\iota
)(a)=a(1\otimes\pi_{g})\right\}  .
\]
where $\iota$\ is the identity automorphism of $B(H_{\pi})$ are, in some
respects more useful$.$ In [14] it is shown that $X_{2}(\pi)$ consists of all
matrices
\[
\left\{  a=\left[  (P_{\alpha})_{ij}(\pi)(x)\right]  (=\left[  a_{ij}\right]
)\in X\otimes B(H_{\pi}):x\in X,1\leq i,j\leq d_{\pi}\right\}  .
\]
It is straightforward to prove that, if $a\in X_{2}(\pi)$ and $x=\sum
_{i}a_{ii}$, then $a_{ij}=(P_{\alpha})_{ij}(\pi)(x).$ In what follows, if
$b\in X\otimes B(H_{\pi})$\ we will denote%
\[
tr(b)=\sum b_{ii}%
\]
which is an $\mathcal{F}$-continuous linear mapping from $X\otimes B(H_{\pi})$
to $X.$ The following lemma is proven for compact non abelian group actions on
C*-algebras in [6, Lemma 2.3.] and for compact abelian $\mathcal{F}$-dynamical
systems in [15]. Since the proof is very similar with the proof of [6, Lemma
2.3.] we will state it without proof

\bigskip

\bigskip\textbf{2.5. Lemma.} \textit{Let (}$X,G,\alpha)$ \textit{be an
}$\mathcal{F}$\textit{-dynamical system with }$G$\textit{\ compact and }%
$J$\textit{\ a two sided ideal of }$X^{\alpha}.$ \textit{Then}%
\[
(\overline{XJX}^{\sigma})^{\alpha}=\mathcal{F}\text{\textit{-closed linear
span of }}\left\{  tr(X_{2}(\pi)JX_{2}(\pi)^{\ast}):\pi\in\widehat{G}\right\}
.
\]
\textit{where, if} $a=\left[  a_{kl}\right]  \in X\otimes B(H_{\pi})$\textit{
and }$j\in X,$\textit{\ by }$ja$\textit{\ we mean the matrix} $\left[
ja_{kl}\right]  $ \textit{and the multiplications }$\overline{XJX}^{\sigma}%
,$\textit{ }$X_{2}(\pi)JX_{2}(\pi)^{\ast}\ $\textit{are defined in 2.6.
below.}

\bigskip

We will use the following notations

\bigskip

\textbf{2.6.} \textbf{Notation.} \textit{Let (}$X,\mathcal{F})$\textit{ be a
dual pair of Banach spaces with }$X$\textit{\ a C*-algebra satisfying
conditions 1)-6). }I\textit{f }$Y,Z$\textit{\ are subsets of }$X$%
\textit{\ denote:}

\textit{a) }$lin\left\{  Y\right\}  $\textit{\ is the linear span of }%
$Y$\textit{.}

\textit{b) }$Y^{\ast}=\left\{  y^{\ast}:y\in Y\right\}  .$

\textit{c) }$YZ=lin\left\{  yz:y\in Y,z\in Z\right\}  .$

\textit{d) }$\overline{Y}^{\sigma}=\mathcal{F}-$\textit{closure of }$Y$ in
$X.$

\textit{e) }$\overline{Y}^{\left\Vert {}\right\Vert }=$\textit{ norm closure
of }$Y.$

\textit{f)} $\overline{Y}^{w}=$ \textit{the }$w^{\ast}$\textit{-closure of
}$Y$\textit{\ in} $\mathcal{F}^{\ast}.$

\textit{If (}$X,G,\alpha)$ \textit{is an }$\mathcal{F}\mathit{-}%
$\textit{dynamical system denote}

\textit{g) }$\mathcal{H}_{\sigma}^{\alpha}(X)$\textit{ the set of all non-zero
globally }$\alpha-$\textit{invariant }$\mathcal{F}\mathit{-}$\textit{closed
hereditary C*}$\mathit{-}$\textit{subalgebras of }$X.$

\bigskip

Notice that if \textit{(}$X,G,\alpha)$ is an $\mathcal{F}-$dynamical system
and if $X_{2}(\pi)$\ is the spectral subspace defined above, then $X_{2}%
(\pi)X_{2}(\pi)^{\ast}$ is a two sided ideal of $X^{\alpha}\otimes B(H_{\pi})$
and $X_{2}(\pi)^{\ast}X_{2}(\pi)\ $is a two sided ideal of $(X\otimes
B(H_{\pi}))^{\alpha\otimes ad\pi}$ where $\alpha\otimes ad\pi$ is the action%
\[
(\alpha_{g}\otimes ad\pi_{g})(a)=(1\otimes\pi_{g})[\alpha_{g}(a_{ij}%
)](1\otimes\pi_{g^{-1}}).
\]
on $X\otimes B(H_{\pi}).$

\bigskip

\textbf{2.7. Definition. }\textit{a)} $sp(\alpha)=\left\{  \pi\in\widehat
{G}:X_{1}(\pi)\neq\left\{  0\right\}  \right\}  .$\newline\textit{b)}
$sp_{\mathcal{F}}(\alpha)=\left\{  \pi\in\widehat{G}:\overline{X_{2}%
(\pi)^{\ast}X_{2}(\pi)}^{\sigma}\text{ \textit{is essential in} }(X\otimes
B(H_{\pi}))^{\alpha\otimes ad\pi}\right\}  .$\newline\textit{c)}
$\widetilde{sp}_{\mathcal{F}}(\alpha)=\left\{  \pi\in\widehat{G}%
:\overline{X_{2}(\pi)^{\ast}X_{2}(\pi)}^{\sigma}=(X\otimes B(H_{\pi}%
))^{\alpha\otimes ad\pi}\right\}  .$\newline\textit{Corresponding to the above
Arveson type spectra b) and c) we define two\ Connes type spectra\newline d)}
$\Gamma_{\mathcal{F}}(\alpha)=\cap\left\{  sp_{\mathcal{F}}(\alpha|_{Y}%
):Y\in\mathcal{H}_{\sigma}^{\alpha}(X)\right\}  .$\newline\textit{e)}
$\widetilde{\Gamma}_{\mathcal{F}}(\alpha)=\cap\left\{  \widetilde
{sp}_{\mathcal{F}}(\alpha|_{Y}):Y\in\mathcal{H}_{\sigma}^{\alpha}(X)\right\}
.$

\bigskip

Clearly, $\widetilde{sp}_{\mathcal{F}}(\alpha)\subset sp_{\mathcal{F}}%
(\alpha)\subset sp(\alpha),$ so $\widetilde{\Gamma}_{\mathcal{F}}%
(\alpha)\subset\Gamma_{\mathcal{F}}(\alpha).$ The definition of $\widetilde
{\Gamma}_{\mathcal{F}}(\alpha)$ is a direct generalization of the strong
Connes spectrum of Kishimoto to compact non abelian groups. Our motivation for
the definition of $\Gamma_{\mathcal{F}}(\alpha)$ above (and $\Gamma(\alpha)$
for C*-dynamical systems in [6]) is the following observation

\bigskip

\textbf{2.8. Remark \textit{a)} }\textit{If (}$X,G,\alpha)$ \textit{is}
\textit{an }$\mathcal{F}\mathit{-}$\textit{dynamical system with }%
$G$\textit{\ compact abelian, then }%
\[
\cap\left\{  sp(\alpha|_{Y}):Y\in\mathcal{H}_{\sigma}^{\alpha}(X)\right\}
=\cap\left\{  sp_{\mathcal{F}}(\alpha|_{Y}):Y\in\mathcal{H}_{\sigma}^{\alpha
}(X)\right\}  .
\]
\textit{and the left hand side of the above equality is the Connes spectrum
for W*}$\mathit{-}$\textit{as well as for C*-dynamical systems.}

\textit{b) If }$G$\textit{\ is not abelian, the equality in part a) is not
true.}

\bigskip

\begin{proof}
a) We have to prove only one inclusion, the opposite one being obvious. Let
$\gamma\in\cap\left\{  sp(\alpha|_{Y}):Y\in\mathcal{H}_{\sigma}^{\alpha
}(X)\right\}  $ and $Y\in\mathcal{H}_{\sigma}^{\alpha}(X).$ Suppose that
$aY_{\gamma}^{\ast}Y_{\gamma}=\left\{  0\right\}  $ for some $a\in Y^{\alpha
},a\neq0.$ Then $aY_{\gamma}^{\ast}=\left\{  0\right\}  .$ Therefore, if we
denote $Z=\overline{aYa^{\ast}}^{\sigma}$, it follows that $Z\in
\mathcal{H}_{\sigma}^{\alpha}(X)$ and $Z_{\gamma}^{\ast}=\left\{  0\right\}  $
which is in contradiction with the hypothesis that $\gamma\in\gamma\in
\cap\left\{  sp(\alpha|_{Y}):Y\in\mathcal{H}_{\sigma}^{\alpha}(X)\right\}
\subset sp(\alpha|_{Z}).$

b) In [14, Example 3.9.] we provided an example of an action of an action of
$G=S_{3}$\ the permutation group on three elements on the algebra $X$ of
$2\times2$ matrices such that $sp(\alpha)=\widehat{G},H_{\sigma}^{\alpha
}=\left\{  X\right\}  ,$ so
\[
\cap\left\{  sp(\alpha|_{Y}):Y\in\mathcal{H}_{\sigma}^{\alpha}(X)\right\}
=sp(\alpha)
\]
\ and we have shown that there exists $\pi\in\widehat{G}$\ such that
($X\otimes B(H_{\pi}))^{\alpha\otimes ad\pi}$ has nontrivial center and,
therefore, it is not a prime C*-algebra. By [6, Thm. 2.2.], it follows that
$\cap\left\{  sp(\alpha|_{Y}):Y\in\mathcal{H}_{\sigma}^{\alpha}(X)\right\}
\neq\cap\left\{  sp_{\mathcal{F}}(\alpha|_{Y}):Y\in\mathcal{H}_{\sigma
}^{\alpha}(X)\right\}  =\Gamma(\alpha).$
\end{proof}

\bigskip

\section{$\mathcal{F}-$ simple fixed point algebras}

\bigskip

Let $(X,G,\alpha)$ be an $\mathcal{F}$-dynamical system with $G$\ compact. In
the rest of this paper we will study how the $\mathcal{F}$-simplicity
(respectively $\mathcal{F}$-primeness) as defined below, of the fixed point
algebras $(X\otimes B(H_{\pi}))^{\alpha\otimes ad\pi}$ is reflected in the
spectral properties of the action.

\bigskip

\textbf{3.1. Definition. }\textit{Let }$(B,\mathcal{F)}$\textit{ be a dual
pair of Banach spaces with }$B$\textit{\ a C*-algebra.}

\textit{a) }$B$\textit{\ is called }$\mathcal{F}$\textit{-simple if every non
zero two sided ideal of }$B$\textit{\ is }$\mathcal{F}$\textit{-dense in }$B.$

\textit{b) }$B$\textit{\ is called }$\mathcal{F}$\textit{-prime if the
annihilator of every non zero two sided ideal of }$B$ \textit{is trivial, or,
equivalently, every non zero two sided ideal of }$B$\textit{\ is an essential
ideal (using Definition 2.1. f) it is easy to see that }$X$\textit{\ is
}$\mathcal{F}-$\textit{prime if and only if }$X$\textit{\ is prime as a
C*-algebra).}

\textit{Let (}$X,G,\alpha)$ \textit{be an }$\mathcal{F}$\textit{-dynamical
system.}

\textit{c) }$X$\textit{\ is called }$\mathcal{\alpha}\mathit{-}$\textit{simple
if every non zero }$\alpha-$\textit{invariant two sided ideal of }%
$X$\textit{\ is }$\mathcal{F}\mathit{-}$\textit{dense in }$X.$

\textit{d) }$X$\textit{\ is called }$\mathcal{\alpha}\mathit{-}$\textit{prime
if every non zero }$\alpha-$\textit{invariant two sided ideal of }%
$X$\textit{\ is an essential ideal. }

\bigskip

In the particular case when $B$\ is a C*-algebra and $\mathcal{F}=B^{\ast}%
$\ is its dual, then, clearly, the concepts of $\mathcal{F}-$simple,
(respectively $\mathcal{F}-$prime) in the above Definition 3.1. a)
(respectively b)) coincide with the usual concepts of simple (respectively
prime) C*-algebras. Similarly, if \textit{(}$X,G,\alpha)$\ is a C*-dynamical
system, that is if $X$\ is a C*-algebra and $\mathcal{F}=X^{\ast}$\ is its
dual, then the notions of $\mathcal{\alpha}\mathit{-}$\textit{simple}\ and
$\mathcal{\alpha}\mathit{-}$\textit{prime}\ coincide with the usual ones for
C*-dynamical systems.

If $B$\ is a von Neumann algebra and $\mathcal{F}=B_{\ast}$\ is its predual,
then, since the weak closure of every essential ideal equals $B,$ it follows
that $B$\ is $\mathcal{F}$-simple if and only if $B$ is $\mathcal{F}-$prime,
so, if and only if $B$\ is a factor. It is also obvious that if \textit{(}%
$X,G,\alpha)\ $is a W*-dynamical system, that is if\ $X$\ is a von Neumann
algebra and $\mathcal{F}=X_{\ast}$ is its predual, then $X$\ is
$\mathcal{\alpha}-$\textit{simpl}e if and only if it is $\mathcal{\alpha}%
-$\textit{prime}, and this holds if and only if $\alpha$\ acts ergodically on
the center of $X$\ (i.e. every fixed element in the center of $X$\ is a scalar).

The above observations and the next Remark show that for W*-dynamical systems,
$(X,G,\alpha)$ with $G$\ compact, the results in the current Section 3 and
Section 4\ are equivalent.

\bigskip

\textbf{3.2. Remark. }\textit{Let (}$X,G,\alpha)$ \textit{be a W*-dynamical
system, that is,\ an }$\mathcal{F}$\textit{-dynamical system with }%
$X$\textit{\ a von Neumann algebra and }$\mathcal{F}=X_{\ast}$\ \textit{its
predual.} \textit{Then} $\widetilde{\Gamma}_{\mathcal{F}}(\alpha
)=\Gamma_{\mathcal{F}}(\alpha).$

\bigskip

\begin{proof}
This follows from the fact that if $X$\ is a von Neumann algebra, $p\in
X^{\alpha}$\ an $\alpha$-invariant projection and $\overline{pX_{2}%
(\pi)p^{\ast}X_{2}(\pi)p}^{\sigma}$ is essential in $(pXp\otimes B(H_{\pi
}))^{\alpha\otimes ad\pi}$, then $\overline{pX_{2}(\pi)p^{\ast}X_{2}(\pi
)p}^{\sigma}$=$(pXp\otimes B(H_{\pi}))^{\alpha\otimes ad\pi}.$
\end{proof}

\bigskip

\bigskip The next lemma will be used in the proofs of the main results of the
current Section 3 and the next Section.

\textbf{3.3. Lemma.} \textit{Let }$(B,G,\alpha)$\textit{\ be an }$\mathcal{F}%
$\textit{-dynamical system with }$G$\textit{\ compact. Then}

\textit{a) If }$\left\{  e_{\lambda}\right\}  $\textit{\ is an approximate
identity of }$B^{\alpha}$\textit{\ in the norm topology, then}%
\[
(norm)\lim_{\lambda}e_{\lambda}x=(norm)\lim_{\lambda}xe_{\lambda}%
=(norm)\lim_{\lambda}e_{\lambda}xe_{\lambda}=x.
\]
\textit{for every }$x\in\overline{\sum_{\pi\in\widehat{G}}B_{1}(\pi
)}^{\left\Vert {}\right\Vert }.$\textit{\ }

\textit{b) If }$b\in B$\textit{\ is such that }$B^{\alpha}bB^{\alpha}=\left\{
0\right\}  $\textit{\ then }$b=0.$

\textit{c) }$\overline{B^{\alpha}BB^{\alpha}}^{\sigma}=\overline{B^{\alpha}%
B}^{\sigma}=\overline{BB^{\alpha}}^{\sigma}=B.$

\textit{d) }$\overline{B^{\alpha}B_{1}(\pi)}^{\sigma}=\overline{B_{1}%
(\pi)B^{\alpha}}^{\sigma}=B_{1}(\pi),\pi\in\widehat{G}.$

\bigskip

\begin{proof}
a) This follows from the proof of [5, Lemma 2.7] in the more general case of
compact quantum group actions.

b) If $\left\{  e_{\lambda}\right\}  $ is an approximate identity of
$B^{\alpha},$\ then $e_{\lambda}be_{\lambda}=0$\ implies%
\[
e_{\lambda}P_{\alpha}(\pi_{ij})(b)e_{\lambda}=P_{\alpha}(e_{\lambda
}be_{\lambda})=0.
\]
for every $\pi\in\widehat{G},1\leq i,j\leq d_{\pi},$ so, by a), $P_{\alpha
}(\pi_{ij})(b)=0$\ Therefore,
\[
\varphi(P_{\alpha}(\pi_{ij})(b))=\int_{G}\pi_{ji}(g)\varphi(\alpha
_{g}(b))dg=0.
\]
for every $\varphi\in\mathcal{F},\pi\in\widehat{G},1\leq i,j\leq d_{\pi}.$
Since $\left\{  \pi_{ij}(g):\pi\in\widehat{G},1\leq i,j\leq d_{\pi}\right\}  $
form an orthogonal basis of $L^{2}(G),$ and $\varphi(\alpha_{g}(b))$ is
continuous on $G,$\ it follows that $\varphi(\alpha_{g}(b))=0$ for every $g\in
G,\varphi\in\mathcal{F},$ so $b=0.$

c) We will prove only that \textit{ }$\overline{B^{\alpha}BB^{\alpha}}%
^{\sigma}=B,$ the proofs of the other equalities being similar.\ Let $\left\{
e_{\lambda}\right\}  $ be an approximate identity of $B^{\alpha}.$ By a),
\[
(norm)\lim_{\lambda}e_{\lambda}xe_{\lambda}=xforeveryx\in\overline{\sum
_{\pi\in\widehat{G}}B_{1}(\pi)}^{\left\Vert {}\right\Vert }.
\]
Therefore%
\[
\overline{\sum_{\pi\in\widehat{G}}B_{1}(\pi)}^{\left\Vert {}\right\Vert
}\subset\overline{B^{\alpha}BB^{\alpha}}^{\left\Vert {}\right\Vert }%
\subset\overline{B^{\alpha}BB^{\alpha}}^{\sigma}.
\]
Since, by Remark 2.4., the $\mathcal{F}$-closure of $\overline{\sum_{\pi
\in\widehat{G}}B_{1}(\pi)}^{\left\Vert {}\right\Vert }$\ equals $B$\ it
follows that
\[
B=\overline{\sum_{\pi\in\widehat{G}}B_{1}(\pi)}^{\sigma}\subset\overline
{B^{\alpha}BB^{\alpha}}^{\sigma}.
\]
so $\overline{B^{\alpha}BB^{\alpha}}^{\sigma}=B.$

d) The proof is similar with the proof of part c).
\end{proof}

\bigskip

Theorem 3.4. below is an extension of [2, Proposition 2.2.2. b)] to the case
of $\mathcal{F}$-dynamical systems with compact groups, not neccessarily
abelian, for the strong Connes spectrum, $\widetilde{\Gamma}_{\mathcal{F}%
}(\alpha).$

\bigskip

\textbf{3.4. Theorem. }\textit{Let }$(X,G,\alpha)$\textit{ be an }%
$F$\textit{-dynamical system with }$G$\textit{\ compact. Then}%
\[
\widetilde{\Gamma}_{\mathcal{F}}(\alpha)=\cap\left\{  \widetilde
{sp}_{\mathcal{F}}(\alpha|_{\overline{JXJ}^{\sigma}}):J\subset X^{\alpha
},\text{ }\mathcal{F}\text{-\textit{closed two sided ideal}}\right\}
\]

\bigskip

\begin{proof}
Clearly, since $\overline{JXJ}^{\sigma}\in\mathcal{H}_{\sigma}^{\alpha}(X),$%
\[
\widetilde{\Gamma}_{\mathcal{F}}(\alpha)\subset\cap\left\{  \widetilde
{sp}_{\mathcal{F}}(\alpha|_{\overline{JXJ}^{\sigma}}):J\subset X^{\alpha
},\text{ }\mathcal{F}\text{-closed two sided ideal}\right\}  .
\]
Let $\pi\in\cap\left\{  \widetilde{sp}_{\mathcal{F}}(\alpha|_{\overline
{JXJ}^{\sigma}}):J\subset X^{\alpha},\text{ }\mathcal{F}\text{-closed two
sided ideal}\right\}  $\ and $Y\in\mathcal{H}_{\sigma}^{\alpha}(X),$ so
$Y^{\alpha}\in\mathcal{H}_{\sigma}(X^{\alpha}).$\ We will prove that
$\overline{Y_{2}(\pi)^{\ast}Y_{2}(\pi)}^{\sigma}=(Y\otimes B(H_{\pi}%
))^{\alpha\otimes ad\pi}$\ and thus $\pi\in\widetilde{sp}_{\mathcal{F}}%
(\alpha|_{Y}).$ Since $Y\in\mathcal{H}_{\sigma}^{\alpha}(X)$\ is arbitrary, it
will follow that $\pi\in\widetilde{\Gamma}_{\mathcal{F}}(\alpha).$ Denote by
$J$\ the following ideal of $X^{\alpha}$%
\[
J=\overline{X^{\alpha}Y^{\alpha}X^{\alpha}}^{\sigma}.
\]
It is clear that $J=\overline{JX^{\alpha}J}^{\sigma}$ (actually it is quite
easy to show that this equality holds without the closure, but we do not need
this fact). Also%
\begin{align}
\overline{Y^{\alpha}JY^{\alpha}}^{\sigma}  &  =\overline{Y^{\alpha}X^{\alpha
}Y^{\alpha}X^{\alpha}Y^{\alpha}}^{\sigma}=\overline{(Y^{\alpha}X^{\alpha
}Y^{\alpha})(Y^{\alpha}X^{\alpha}Y^{\alpha})}^{\sigma}=\label{1}\\
&  =\overline{Y^{\alpha}Y^{\alpha}}^{\sigma}=Y^{\alpha}.\nonumber
\end{align}
Denote $Z=\overline{JXJ}^{\sigma}.$ Notice that, since $Y\in\mathcal{H}%
_{\sigma}^{\alpha}(X),$ we have $\overline{Y^{\alpha}XY^{\alpha}}^{\sigma}=Y,$
so
\begin{gather}
Z=\overline{X^{\alpha}Y^{\alpha}X^{\alpha}}^{\sigma}X\overline{X^{\alpha
}Y^{\alpha}X^{\alpha}}^{\sigma}=\overline{X^{\alpha}Y^{\alpha}X^{\alpha
}XX^{\alpha}Y^{\alpha}X^{\alpha}}^{\sigma}=\label{2}\\
\overline{X^{\alpha}Y^{\alpha}(X^{\alpha}XX^{\alpha})Y^{\alpha}X^{\alpha}%
}^{\sigma}=\overline{X^{\alpha}Y^{\alpha}XY^{\alpha}X^{\alpha}}^{\sigma
}=\overline{X^{\alpha}YX^{\alpha}}^{\sigma}.\nonumber
\end{gather}
Since $\pi\in\cap\left\{  \widetilde{sp}_{\mathcal{F}}(\alpha|_{\overline
{JXJ}^{\sigma}}):J\subset X^{\alpha},\text{ }\mathcal{F}\text{-closed two
sided ideal}\right\}  $, it follows that $\pi\in\widetilde{sp}_{\mathcal{F}%
}(\alpha|_{Z}),$ so
\begin{equation}
\overline{Z_{2}(\pi)^{\ast}Z_{2}(\pi)}^{\sigma}=(Z\otimes B(H_{\pi}%
))^{\alpha\otimes ad\pi}. \label{3}%
\end{equation}
Using the equalities (\ref{2}) above, the fact that $Y$\ is a hereditary
C*-subalgebra of $X$, and the obvious equality%
\[
P_{ij}(\pi)(xyz)=xP_{ij}(\pi)(y)z
\]
for every $x,z\in X^{\alpha},y\in X,$\ and $1\leq i,j\leq\dim H_{\pi},$\ the
relation (\ref{3}) becomes%
\begin{equation}
\overline{X^{\alpha}Y_{2}(\pi)^{\ast}Y_{2}(\pi)X^{\alpha}}^{\sigma}%
=\overline{X^{\alpha}(Y\otimes B(H_{\pi}))^{\alpha\otimes ad\pi}X^{\alpha}%
}^{\sigma}. \label{4}%
\end{equation}
where, for $x\in X^{\alpha}$\ and $a\in X\otimes B(H_{\pi}),a=\left[
a_{kl}\right]  ,$\ by $xa$ we mean the matrix whose $kl$\ entry is $xa_{kl}.$
Therefore, by applying Lemma 3.3. d) to $B=Y,$\ we get%
\begin{equation}
\overline{X^{\alpha}Y^{\alpha}Y_{2}(\pi)^{\ast}Y_{2}(\pi)Y^{\alpha}X^{\alpha}%
}^{\sigma}=\overline{X^{\alpha}Y^{\alpha}(Y\otimes B(H_{\pi}))^{\alpha\otimes
ad\pi}Y^{\alpha}X^{\alpha}}^{\sigma}. \label{5}%
\end{equation}
By multiplying (\ref{5}) on the right and on the left by $Y^{\alpha}$\ and
taking into account that, by Lemma 3.3. c) $\overline{Y^{\alpha}YY^{\alpha}%
}^{\sigma}=Y$ and consequently, $\overline{Y^{\alpha}X^{\alpha}Y^{\alpha}%
}^{\sigma}=Y^{\alpha},$ it follows that
\[
\overline{Y_{2}(\pi)^{\ast}Y_{2}(\pi)}^{\sigma}=(Y\otimes B(H_{\pi}%
))^{\alpha\otimes ad\pi}.
\]
Therefore, $\pi\in\widetilde{sp}_{\mathcal{F}}(\alpha|_{Y})$ and the proof is complete.
\end{proof}

\bigskip

In the next Lemma and the rest of the paper, a subalgebra of $X\otimes
B(H_{\pi})$\ will be called $\mathcal{F}-$simple\ (respectively $\mathcal{F}%
-$prime) if it is $\mathcal{F}\otimes B(H_{\pi})^{\ast}-$simple (respectively
$\mathcal{F}\otimes B(H_{\pi})^{\ast}-$prime) where $B(H_{\pi})^{\ast}$
denotes the dual of $B(H_{\pi}).$ Clearly, a subalgebra of $X\otimes B(H_{\pi
})$ is $\mathcal{F}-$prime if and only if it is a prime C*-algebra. The
similar statement for the $\mathcal{F}-$simple case is not true.

\bigskip

\textbf{3.5. Lemma. }\textit{Let }$(X,G,\alpha)$\textit{ be an }$\mathcal{F}%
$\textit{-dynamical system with }$G$\textit{\ compact. Then, if }$X^{\alpha}%
$\textit{\ is }$\mathcal{F}$\textit{-simple, it follows that (}$X\otimes
B(H_{\pi}))^{\alpha\otimes ad\pi}$\textit{\ is }$\mathcal{F}-$\textit{simple.
}

\bigskip

\begin{proof}
Let $\pi\in\widetilde{sp}_{\mathcal{F}}(\alpha)\subset sp(\alpha).$\ Since
$X^{\alpha}$\ is $\mathcal{F}$-simple, so $X^{\alpha}\otimes B(H_{\pi})$ is
also $\mathcal{F}$-simple and $X_{2}(\pi)X_{2}(\pi)^{\ast}$ is an ideal of
$X^{\alpha}\otimes B(H_{\pi}),$\ it follows that $\overline{X_{2}(\pi
)X_{2}(\pi)^{\ast}}^{\sigma}=X^{\alpha}\otimes B(H_{\pi}).$ To prove that
($X\otimes B(H_{\pi}))^{\alpha\otimes ad\pi}$\ is simple, let $I\subset
(X\otimes B(H_{\pi}))^{\alpha\otimes ad\pi}$ be a non-zero ideal. Then it can
be easily verified that%
\begin{align*}
J  &  =\overline{lin}^{\sigma}\left\{  yy^{\ast}:y\in X_{2}(\pi)I\right\}  =\\
&  =\overline{X_{2}(\pi)IX_{2}(\pi)^{\ast}}^{\sigma}.
\end{align*}
is an ideal of $X^{\alpha}\otimes B(H_{\pi})$\ and, since the latter algebra
is $\mathcal{F}$-simple, it follows that $J=X^{\alpha}\otimes B(H_{\pi}).$
Therefore, since $\pi\in\widetilde{sp}_{\mathcal{F}}(\alpha),$ we have
$\overline{X_{2}(\pi)^{\ast}X_{2}(\pi)}^{\sigma}=(X\otimes B(H_{\pi}%
))^{\alpha\otimes ad\pi}$ and consequently, since, by Lemma 3.3. d)
$\overline{X^{\alpha}X_{2}(\pi)}^{\sigma}=X_{2}(\pi),$ we have
\[
(X\otimes B(H_{\pi}))^{\alpha\otimes ad\pi}=\overline{X_{2}(\pi)^{\ast}%
JX_{2}(\pi)}^{\sigma}\subset\overline{X_{2}(\pi)^{\ast}X_{2}(\pi)IX_{2}%
(\pi)^{\ast}X_{2}(\pi)}^{\sigma}\subset I.
\]
Thus $I=(X\otimes B(H_{\pi}))^{\alpha\otimes ad\pi}$ and we are done.
\end{proof}

\bigskip

The following result extends [2, Th\'{e}or\`{e}me 2.4.1], [12, Theorem 2.
i)$\Leftrightarrow ii)]$] and [15, Theorem 3.4.] to the more general case of
$\mathcal{F}$-dynamical systems and non abelian compact groups $G.$

\bigskip

\textbf{3.6. Theorem. }\textit{Let }$(X,G,\alpha)$\textit{ be an }%
$\mathcal{F}$\textit{-dynamical system with }$G$\textit{\ compact. The
following conditions are equivalent:\newline i) (}$X\otimes B(H_{\pi
}))^{\alpha\otimes ad\pi}$\textit{ is }$\mathcal{F}$\textit{-simple for all
}$\pi\in sp(\alpha).$\textit{\newline ii) }$X$\textit{ is }$\alpha
$\textit{-simple and }$sp(\alpha)=\widetilde{\Gamma}_{\mathcal{F}}(\alpha).$

\bigskip

\begin{proof}
$i)\Rightarrow ii)$ Suppose that ($X\otimes B(H_{\pi}))^{\alpha\otimes ad\pi}$
is $\mathcal{F}$-simple for all $\pi\in sp(\alpha)$. Then, it follows
immediately from the definitions that $sp(\alpha)=\widetilde{sp}_{\mathcal{F}%
}(\alpha).$ Let $\pi\in sp(\alpha)$ be arbitrary. Since, in particular,
$X^{\alpha}$\ is $\mathcal{F}$-simple, so it\ has no non-trivial $\mathcal{F}%
$-closed ideals, from Theorem 3.4.\ it follows that $\pi\in\widetilde{\Gamma
}_{\mathcal{F}}(\alpha),$ so $sp(\alpha)=\widetilde{\Gamma}_{\mathcal{F}%
}(\alpha).$ Let us prove that $X$ \ is\ $\alpha$-simple.\ If $I$\ is an
$\mathcal{F}$-closed $\alpha$-invariant ideal of $X,$\ then $I^{\alpha}$\ is
an $\mathcal{F}$-closed ideal of $X^{\alpha}$, so $I^{\alpha}=X^{\alpha}.$ By
Lemma 3.3. c) applied to $B=I,$ and to $B=X$\ it follows that $\overline
{I^{\alpha}II^{\alpha}}^{\sigma}=I$ and $\overline{X^{\alpha}XX^{\alpha}%
}^{\sigma}=X$,\ so%
\[
X=\overline{X^{\alpha}XX^{\alpha}}^{\sigma}=\overline{I^{\alpha}XI^{\alpha}%
}^{\sigma}\subset I.
\]
Therefore, $\ I=X,$ hence $X$ is $\alpha$-simple.\newline$ii)\Rightarrow i).$
Suppose that $X$\ is $\alpha$-simple and $sp(\alpha)=\widetilde{\Gamma
}_{\mathcal{F}}(\alpha).$ We will prove first that $X^{\alpha}$\ is
$\mathcal{F}$-simple. Let $J\subset X^{\alpha}$ be a non zero ideal and
$\pi\in\widetilde{\Gamma}_{\mathcal{F}}(\alpha).$ Since $\overline
{JXJ\ }^{\sigma}\in\mathcal{H}_{\sigma}^{\alpha}(X),$ and $\pi\in
\widetilde{\Gamma}_{\mathcal{F}}(\alpha),$ it follows that
\begin{equation}
\overline{JX_{2}(\pi)^{\ast}JX_{2}(\pi)J}^{\sigma}=\overline{J(X\otimes
B(H_{\pi}))^{\alpha\otimes ad\pi}J}^{\sigma}. \label{6}%
\end{equation}
where, for $j\in J\subset X^{\alpha}$\ and $a\in X\otimes B(H_{\pi}),a=\left[
a_{kl}\right]  ,$\ by $ja$ we mean the matrix whose $kl$\ entry is $ja_{kl}%
.$\ By multiplying the above relation on the left by $X_{2}(\pi)$\ and on the
right by $X_{2}(\pi)^{\ast},$\ we get%
\begin{equation}
\overline{X_{2}(\pi)JX_{2}(\pi)^{\ast}JX_{2}(\pi)JX_{2}(\pi)^{\ast}}^{\sigma
}=\overline{X_{2}(\pi)J(X\otimes B(H_{\pi}))^{\alpha\otimes ad\pi}JX_{2}%
(\pi)^{\ast}}^{\sigma}. \label{7}%
\end{equation}
From the above relations (\ref{6}) and (\ref{7}) it follows that%
\[
X_{2}(\pi)JX_{2}(\pi)^{\ast}=\overline{X_{2}(\pi)JX_{2}(\pi)^{\ast}}^{\sigma
}\subset\overline{X_{2}(\pi)J(X\otimes B(H_{\pi}))^{\alpha\otimes ad\pi}%
JX_{2}(\pi)^{\ast}}^{\sigma}=
\]%
\[
=\overline{X_{2}(\pi)JX_{2}(\pi)^{\ast}JX_{2}(\pi)JX_{2}(\pi)^{\ast}}^{\sigma
}\subset\overline{(X^{\alpha}\otimes B(H_{\pi}))J(X^{\alpha}\otimes B(H_{\pi
}))}^{\sigma}\subset
\]%
\[
J\otimes B(H_{\pi})
\]
It follows that $tr(X_{2}(\pi)JX_{2}(\pi)^{\ast})\subset J.$ From Lemma 2.5.
it follows that $(\overline{XJX}^{\sigma})^{\alpha}\subset J.$ Since $X$\ is
$\alpha$-simple, we have $\overline{XJX}^{\sigma}=X,$ so $J=X^{\alpha}$ and
therefore, $X^{\alpha}$ is $\mathcal{F}$-simple. Applying lemma 3.5. it
follows that ($X\otimes B(H_{\pi}))^{\alpha\otimes ad\pi}$\ is $\mathcal{F}%
$-simple for all $\pi\in sp(\alpha)=\widetilde{\Gamma}_{\mathcal{F}}(\alpha).$
\end{proof}

\bigskip

\section{$\mathcal{F}$-prime fixed point algebras}

\bigskip

This section is concerned with the relationship between the $\mathcal{F}%
$-primeness of the fixed point algebras and the spectral properties, involving
the Connes spectrum $\Gamma_{\mathcal{F}}(\alpha)$\ of the $\mathcal{F}%
$-dynamical system $(X,G,\alpha).$

\bigskip

Theorem 4.1. below is an extension of [2, Proposition 2.2.2. b)] to the case
of $\mathcal{F}$-dynamical systems with compact groups, not neccessarily
abelian, for the Connes spectrum, $\Gamma_{\mathcal{F}}(\alpha).$ By Remark
3.2. and the discussion preceding it, if $(X,G,\alpha)$ is a W*-dynamical
system (that is $X$\ is a von Neumann algebra and $\mathcal{F}=X_{\ast}$ its
predual), then the next Theorem 4.1. is equivalent with Theorem 3.4.

\bigskip

\textbf{4.1. Theorem}. \textit{Let }$(X,G,\alpha)$\textit{ be an}
$\mathcal{F}$-\textit{dynamical system}. \textit{Then}%
\[
\Gamma_{\mathcal{F}}(\alpha)=\cap\left\{  sp_{\mathcal{F}}(\alpha
|_{\overline{JXJ}^{\sigma}}):J\subset X^{\alpha},\text{ \textit{a non-zero}
}\mathcal{F}\text{-\textit{closed two sided ideal}}\right\}  \text{.}%
\]

\bigskip

\begin{proof}
Since $\overline{JXJ}^{\sigma}\in\mathcal{H}_{\sigma}^{\alpha}(X)$ for every
non-zero $\mathcal{F}$-closed two sided ideal $J\subset X^{\alpha}$, we
have$,$%
\[
\Gamma_{\mathcal{F}}(\alpha)\subset\cap\left\{  sp_{_{\mathcal{F}}}%
(\alpha|_{\overline{JXJ}^{\sigma}}):J\subset X^{\alpha},\text{ }%
\mathcal{F}\text{-closed two sided ideal}\right\}  .
\]
Now let $\pi\in\cap\left\{  sp_{_{\mathcal{F}}}(\alpha|_{\overline
{JXJ}^{\sigma}}):J\subset X^{\alpha},\text{ }\mathcal{F}\text{-closed two
sided ideal}\right\}  $ and $Y\in\mathcal{H}_{\sigma}^{\alpha}(X),$ so
$Y^{\alpha}\in\mathcal{H}_{\sigma}(X^{\alpha}).$\ We will prove that
$\overline{Y_{2}(\pi)^{\ast}Y_{2}(\pi)}^{\sigma}$\ is essential in $(Y\otimes
B(H_{\pi}))^{\alpha\otimes ad\pi}.$ As in the proof of Theorem 3.4., let
$J=\overline{X^{\alpha}Y^{\alpha}X^{\alpha}}^{\sigma}$\ and $Z=\overline
{JXJ}^{\sigma}\in H_{\sigma}^{\alpha}(X).$\ Since $\pi\in\cap\left\{
sp_{_{\mathcal{F}}}(\alpha|_{\overline{JXJ}^{\sigma}}):J\subset X^{\alpha
},\text{ }\mathcal{F}\text{-closed two sided ideal}\right\}  ,$\ we have that
\ $\pi\in sp_{_{\mathcal{F}}}(\alpha|_{Z}).$ Therefore, $\overline{Z_{2}%
(\pi)^{\ast}Z_{2}(\pi)}^{\sigma}$\ is essential in $(Z\otimes B(H_{\pi
}))^{\alpha\otimes ad\pi}.$\ As noticed in the proof of Theorem 3.4.,
\[
\overline{Z_{2}(\pi)^{\ast}Z_{2}(\pi)}^{\sigma}=\overline{X^{\alpha}Y_{2}%
(\pi)^{\ast}Y_{2}(\pi)X^{\alpha}}^{\sigma}.
\]
and
\[
(Z\otimes B(H_{\pi}))^{\alpha\otimes ad\pi}=\overline{X^{\alpha}(Y\otimes
B(H_{\pi}))^{\alpha\otimes ad\pi}X^{\alpha}}^{\sigma}\text{.}%
\]
Let $a\in(Y\otimes B(H_{\pi}))^{\alpha\otimes ad\pi}$\ be such that
\[
a\overline{Y_{2}(\pi)^{\ast}Y_{2}(\pi)}^{\sigma}=\left\{  0\right\}  .
\]
Then, by Lemma 3.3. c) $Y=\overline{Y^{\alpha}Y}^{\sigma}=\overline{Y^{\alpha
}X^{\alpha}Y^{\alpha}Y}^{\sigma}$ and, by Lemma 3.3. d), $Y_{2}(\pi)^{\ast
}=\overline{Y^{\alpha}Y_{2}(\pi)^{\ast}}^{\sigma}.$\ Then, it follows that%
\[
a\overline{Y_{2}(\pi)^{\ast}Y_{2}(\pi)}^{\sigma}=\overline{aY^{\alpha}%
Y_{2}(\pi)^{\ast}Y_{2}(\pi)}^{\sigma}=\overline{aY^{\alpha}X^{\alpha}%
Y^{\alpha}Y_{2}(\pi)^{\ast}Y_{2}(\pi)}^{\sigma}=\left\{  0\right\}  .
\]
Therefore%
\[
\overline{Y^{\alpha}aY^{\alpha}X^{\alpha}Y^{\alpha}Y_{2}(\pi)^{\ast}Y_{2}%
(\pi)X^{\alpha}}^{\sigma}=\left\{  0\right\}  .
\]
so, since $\overline{Z_{2}(\pi)^{\ast}Z_{2}(\pi)}^{\sigma}=\overline
{X^{\alpha}Y_{2}(\pi)^{\ast}Y_{2}(\pi)X^{\alpha}}^{\sigma}$\ is essential in
$\overline{X^{\alpha}(Y\otimes B(H_{\pi}))^{\alpha\otimes ad\pi}X^{\alpha}%
}^{\sigma},$ we have $Y^{\alpha}aY^{\alpha}=0$ and therefore, by Lemma 3.3. b)
applied to $B=Y,$ it follows that $a=0.$
\end{proof}

\bigskip

\textbf{4.2. Lemma. }\textit{Let }$(X,G,\alpha)$\textit{ be an }$\mathcal{F}%
$\textit{-dynamical system with }$G$\textit{\ compact. Then, if }$X^{\alpha}%
$\textit{\ is }$\mathcal{F}$\textit{-prime, it follows that (}$X\otimes
B(H_{\pi}))^{\alpha\otimes ad\pi}$\textit{\ is }$\mathcal{F}-$\textit{prime
for every }$\pi\in sp_{\mathcal{F}}(\alpha)$\textit{.}

\bigskip

\begin{proof}
Since $X^{\alpha}$\ is $\mathcal{F}$-prime. it follows that $X^{\alpha}\otimes
B(H_{\pi})$ is $\mathcal{F}$-prime for every $\pi\in\widehat{G}$. Since
$X_{2}(\pi)X_{2}(\pi)^{\ast}$ is a non zero ideal of $X^{\alpha}\otimes
B(H_{\pi}),$ it follows that $X_{2}(\pi)X_{2}(\pi)^{\ast}$is an essential
ideal$.$ To prove that ($X\otimes B(H_{\pi}))^{\alpha\otimes ad\pi}$\ is
$\mathcal{F}$-prime, let $I\subset(X\otimes B(H_{\pi}))^{\alpha\otimes ad\pi}$
be a non-zero ideal. Then, as in the proof of Lemma 3.5., consider the
following ideal of $X^{\alpha}\otimes B(H_{\pi})$
\begin{align*}
J  &  =\overline{lin}^{\sigma}\left\{  yy^{\ast}:y\in X_{2}(\pi)I\right\}  =\\
&  =\overline{X_{2}(\pi)IX_{2}(\pi)^{\ast}}^{\sigma}.
\end{align*}
Since $X^{\alpha}\otimes B(H_{\pi})$ is $\mathcal{F}$-prime, it follows that
$J$ is essential in $X^{\alpha}\otimes B(H_{\pi}).$ Therefore, if
$a\in(X\otimes B(H_{\pi}))^{\alpha\otimes ad\pi}$ and $aI=\left\{  0\right\}
.$\ we have
\[
aX_{2}(\pi)^{\ast}JX_{2}(\pi)\subset\overline{aX_{2}(\pi)^{\ast}X_{2}%
(\pi)IX_{2}(\pi)^{\ast}X_{2}(\pi)}^{\sigma}\subset\overline{aI}^{\sigma
}=\left\{  0\right\}  .
\]
so
\[
(X_{2}(\pi)aX_{2}(\pi)^{\ast}J)X_{2}(\pi)X_{2}(\pi)^{\ast}=\left\{  0\right\}
\]
Thus, since $X_{2}(\pi)X_{2}(\pi)^{\ast}$ is essential in $X^{\alpha}\otimes
B(H_{\pi}),\ $it\ follows that $X_{2}(\pi)aX_{2}(\pi)^{\ast}J=\left\{
0\right\}  .$ Since $X_{2}(\pi)aX_{2}(\pi)^{\ast}\subset X^{\alpha}\otimes
B(H_{\pi}),$ $J$\ is essential in $X^{\alpha}\otimes B(H_{\pi})$ and $\pi\in
sp_{\mathcal{F}}(\alpha),$ it follows that $a=0.$
\end{proof}

\bigskip

The next result extends [2, Th\'{e}or\`{e}me 2.4.1.], and [13, Theorem
8.10.4.] to the case of $\mathcal{F}$-dynamical systems with compact non
abelian groups.

\bigskip\ 

\textbf{4.3. Theorem. }\textit{Let }$(X,G,\alpha)$\textit{ be an }%
$\mathcal{F}$\textit{-dynamical system with }$G$\textit{\ compact. The
following conditions are equivalent:\newline i) (}$X\otimes B(H_{\pi
}))^{\alpha\otimes ad\pi}$\textit{ is }$\mathcal{F}$\textit{-prime for all
}$\pi\in sp(\alpha).$\textit{\newline ii) }$X$\textit{ is }$\alpha
$\textit{-prime and }$sp(\alpha)=\Gamma_{\mathcal{F}}(\alpha).$

\bigskip

\begin{proof}
$i)\Rightarrow ii)$ Suppose that ($X\otimes B(H_{\pi}))^{\alpha\otimes ad\pi}$
is $\mathcal{F}$-prime for all $\pi\in sp(\alpha)$. Then, it follows
immediately from i) and the definitions that $sp(\alpha)=sp_{\mathcal{F}%
}(\alpha).$ Let $\pi\in sp(\alpha)$ be arbitrary. We will use Theorem 4.1. to
show that $\pi\in\Gamma_{\mathcal{F}}(\alpha)$. Indeed, let $J$ be a
non-trivial ideal of $X^{\alpha}$ and $Z=\overline{JXJ}^{\sigma}\in H_{\sigma
}^{\alpha}(X).$\ We will show that $\pi\in sp(\alpha|_{Z})$, that is
$Z_{2}(\pi)^{\ast}Z_{2}(\pi)\ $is.esential in ($Z\otimes B(H_{\pi}%
))^{\alpha\otimes ad\pi}.$ Notice that%
\[
Z_{2}(\pi)=JX_{2}(\pi)J.
\]
so%
\[
Z_{2}(\pi)^{\ast}Z_{2}(\pi)=JX_{2}(\pi)^{\ast}JX_{2}(\pi)J.
\]
Since, in particular, $X^{\alpha}$ is prime, and $J\ $is an essential ideal of
$X^{\alpha},$ we have$\ Z_{2}(\pi)\neq\left\{  0\right\}  .$ Indeed as
observed after 2.6.$,X_{2}(\pi)X_{2}(\pi)^{\ast}$ is an ideal of $X^{\alpha
}\otimes B(H_{\pi}),$ so, as $X^{\alpha}$\ is $\mathcal{F}$-prime\ it follows
that $X^{\alpha}\otimes B(H_{\pi})$ is a prime C*-algebra$\ $\ and
therefore\ $JX_{2}(\pi)X_{2}(\pi)^{\ast}\neq\left\{  0\right\}  ,$ hence
$JX_{2}(\pi)\neq\left\{  0\right\}  $ and $X_{2}(\pi)^{\ast}J\neq\left\{
0\right\}  ,$\ so, $X_{2}(\pi)^{\ast}JX_{2}(\pi)\neq\left\{  0\right\}
.$\ Using the hypothesis that \textit{(}$X\otimes B(H_{\pi}))^{\alpha\otimes
ad\pi}$\ is $\mathcal{F}-$prime and the fact that $J$\ is a non-trivial ideal
of $X^{\alpha},$\ it follows that $X_{2}(\pi)^{\ast}JX_{2}(\pi)J\neq\left\{
0\right\}  ,$ so, $Z_{2}(\pi)\neq\left\{  0\right\}  .$\ As noticed\ above,
$X_{2}(\pi)^{\ast}JX_{2}(\pi)$\ is a non-trivial ideal of ($X\otimes B(H_{\pi
}))^{\alpha\otimes ad\pi}.$\ If $a\in$($Z\otimes B(H_{\pi}))^{\alpha\otimes
ad\pi},a\geq0$\ is such that $aZ_{2}(\pi)^{\ast}Z_{2}(\pi)=\left\{  0\right\}
,$ then%
\[
aJX_{2}(\pi)^{\ast}JX_{2}(\pi)J=\left\{  0\right\}  .
\]
Hence%
\[
JaJX_{2}(\pi)^{\ast}JX_{2}(\pi)J=\left\{  0\right\}  .
\]
Since, as noticed above, $X_{2}(\pi)^{\ast}JX_{2}(\pi)$\ is non-trivial and
($X\otimes B(H_{\pi}))^{\alpha\otimes ad\pi}$ is $\mathcal{F}$-prime it
follows that
\[
JaJ=\left\{  0\right\}  .
\]
so $Ja=\left\{  0\right\}  .$ Hence $Jtr(a)=\left\{  0\right\}  .$\ Since
$X^{\alpha}$ is $\mathcal{F}$-prime,\ we deduce that $tr(a)=0,$ so
$a=0$\ because\ $a$ was assumed to be non negative. Therefore, $\pi\in
\Gamma_{\mathcal{F}}(\alpha),$ so $sp(\alpha)=\Gamma_{\mathcal{F}}(\alpha).$
It remains to prove that $X$\ is $\mathcal{F}$-prime. Let $I\subset X$\ be an
$\alpha$-invariant\ non-trivial ideal and $x\in X,x\geq0$ be such that
$xI=\left\{  0\right\}  .$ Then, in particular, $xI^{\alpha}=\left\{
0\right\}  ,$\ so $P(x)I^{\alpha}=\left\{  0\right\}  .$ Since $X^{\alpha}%
$\ is $\mathcal{F}$-prime and $I^{\alpha}$ is a non trivial ideal of
$X^{\alpha}$\ we have $P(x)=0$\ so, since $P\ $is faithful, $x=0.$\newline
ii)$\Rightarrow i)$ Suppose that $X$\ is $\alpha$-prime and $sp(\alpha
)=\Gamma_{\mathcal{F}}(\alpha).$\ We will prove first that $X^{\alpha}$\ is
$\mathcal{F}$-prime. Let $J\subset X^{\alpha}$\ be a non zero ideal and $a\in
X^{\alpha},a\geq0,a\neq0$\ such that $Ja=\left\{  0\right\}  .$ Since $X$\ is
$\alpha-$prime, and $XJX$\ is a non zero $\alpha$-invariant ideal of $X,$ it
follows that $XJXa\neq\left\{  0\right\}  $ so $JXa\neq\left\{  0\right\}
.$\ Therefore, since by Remark 2.4.,
\[
X=\overline{\sum_{\pi\in sp(\alpha)}X_{1}(\pi)}^{\sigma}.
\]
there exists $\pi\in sp(\alpha)$\ such that $JX_{1}(\pi)a\neq\left\{
0\right\}  .$ Denote $Z=\overline{aXa}^{\sigma}\in H_{\sigma}^{\alpha}(X).$
Then, since $\pi\in sp(\alpha)=\Gamma_{\mathcal{F}}(\alpha),$ $Z_{2}%
(\pi)^{\ast}Z_{2}(\pi)$ is essential in ($Z\otimes B(H_{\pi}))^{\alpha\otimes
ad\pi}.$ But%
\[
\overline{Z_{2}(\pi)^{\ast}Z_{2}(\pi)}^{\sigma}=\overline{aX_{2}(\pi)^{\ast
}a^{2}X_{2}(\pi)a}^{\sigma}.
\]
and%
\[
(Z\otimes B(H_{\pi}))^{\alpha\otimes ad\pi}=\overline{a((X\otimes B(H_{\pi
}))^{\alpha\otimes ad\pi})a}^{\sigma}.
\]
Taking into account that $X_{2}(\pi)a^{2}X_{2}(\pi)^{\ast}\subset X^{\alpha
}\otimes B(H_{\pi})$\ and $Ja=\left\{  0\right\}  $ we immediately get that%

\[
JX_{2}(\pi)a^{2}X_{2}(\pi)^{\ast}\subset J\otimes B(H_{\pi}).
\]
Hence%

\[
\overline{(aX_{2}(\pi)^{\ast}J)X_{2}(\pi)a^{2}X_{2}(\pi)^{\ast}(a^{2}X_{2}%
(\pi)a)}^{\sigma}=\left\{  0\right\}  .
\]
Therefore%

\[
\overline{(aX_{2}(\pi)^{\ast}JX_{2}(\pi)a)(aX_{2}(\pi)^{\ast}a^{2}X_{2}%
(\pi)a)}^{\sigma}=\left\{  0\right\}  .
\]
It follows that%
\[
(aX_{2}(\pi)^{\ast}JX_{2}(\pi)a)(Z_{2}(\pi)^{\ast}Z_{2}(\pi))=\left\{
0\right\}  .
\]
Since $Z_{2}(\pi)^{\ast}Z_{2}(\pi)$\ is essential in $(Z\otimes B(H_{\pi
}))^{\alpha\otimes ad\pi}$\ and obviously $aX_{2}(\pi)^{\ast}JX_{2}(\pi
)^{\ast}a\subset(Z\otimes B(H_{\pi}))^{\alpha\otimes ad\pi}$ it follows that
$JX_{2}(\pi)a=\left\{  0\right\}  $ and hence $JX_{1}(\pi)a=\left\{
0\right\}  ,$ but this is in contradiction with our choice of $\pi$ in
$sp(\alpha)$ so $X^{\alpha}$\ is $\mathcal{F}$-prime. From Lemma 4.2. it
follows that ($X\otimes B(H_{\pi}))^{\alpha\otimes ad\pi}$ is $\mathcal{F}%
$-prime for all $\pi\in\Gamma_{\mathcal{F}}(\alpha)=sp_{\mathcal{F}}%
(\alpha)=sp(\alpha)$ and we are done.
\end{proof}

\bigskip

\bigskip

{\LARGE References}

\bigskip

[1] W.B. Arveson, On groups of automorphisms of operator algebras, J. Funct.
Anal. 15 (1974) 217-243.

[2] A. Connes, Une classification des facteurs de type III, Ann. Sci. \'{E}c.
Norm. Sup\'{e}r. 6 (1973) 133-252.

[3] A. Connes, M. Takesaki, The flow of weights on factors of type III,
T\^{o}hoku Math. J. 29 (1977) 473-575.

[4] C. D'Antoni, L. Zsido, Groups of linear isometries on multiplier
C*-algebras, Pacific J. Math. 193 (2000) 279-306.

[5] R. Dumitru, C. Peligrad, Compact quantum group actions and invariant
derivations, Proc. Amer. Math. Soc. 135 (2007) 3977-3984.

[6] E. C. Gootman, A. J. Lazar, C. Peligrad, Spectra for compact group
actions, J. Operator Theory 31 (1994) 381-399.

[7] E. Hille, R. Phillips, Functional Analysis and Semi-Groups, AMS, 1957.

[8] A. Kishimoto, Simple crossed products of C*-algebras by locally compact
abelian groups, Yokohama Math. J. 28 (1980) 69-85.

[9] M. B. Landstad, Algebras of spherical functions associated with covariant
systems over a compact group, Math. Scand. 47 (1980), 137-149.

[10] D. Olesen, Inner $^{\ast}$-automorphisms of simple C*-algebras, Comm.
Math. Phys..44 (1975) 175-190.

[11] D. Olesen, G.K. Pedersen, Applications of the Connes spectrum to
C*-dynamical systems, J. Funct. Anal. 30 (1978) 179-197.

[12] D. Olesen, G.K. Pedersen, E. Stormer, Compact abelian groups of
automorphisms of simple C*-algebras, Invent. Math. 39 (1977) 55-64.

[13] G.K. Pedersen, C*-algebras and Their Automorphism Groups, Academic Press, 1979.

[14] C. Peligrad, Locally compact group actions and compact subgroups, J.
Funct. Anal. 76 (1988) 126-139.

[15] C. Peligrad, A solution of the maximality problem for one-parameter
dynamical systems, Adv. Math. 329 (2018) 742-780.

[16] L. Zsido, Spectral and ergodic properties of the analytic generators, J.
Approx. Theory 20 (1977) 77-138.

\end{document}